\documentclass{amsart}

\usepackage{amssymb}

\newtheorem{thm}{Theorem}

\newtheorem{lem}[thm]{Lemma}
\newtheorem{cor}[thm]{Corollary}

\newtheorem{prop}[thm]{Proposition}

\theoremstyle{definition}
\newtheorem{defn}[thm]{Definition}

\newtheorem{say}[thm]{}
\newtheorem{exmp}[thm]{Example}

\newtheorem{rem}[thm]{Remark}          

\newtheorem{ack}{Acknowledgments}

\newtheorem{defn-thm}[thm]{Definition--Theorem}  
\newtheorem{defn-lem}[thm]{Definition--Lemma}  

\theoremstyle{remark}

\renewcommand{\o}[0]{{\mathcal O}} 

\renewcommand{\a}[0]{{\mathbb A}}

\newcommand{\p}[0]{{\mathbb P}}

\newcommand{\qtq}[1]{\quad\mbox{#1}\quad}
\newcommand{\spec}[0]{\operatorname{Spec}}

\newcommand{\supp}[0]{\operatorname{Supp}}    
\newcommand{\red}[0]{\operatorname{red}}    
\newcommand{\codim}[0]{\operatorname{codim}}    
\newcommand{\im}[0]{\operatorname{im}}

\newcommand{\coker}[0]{\operatorname{coker}}    
    
\newcommand{\Hom}[0]{\operatorname{Hom}}

\newcommand{\tors}[0]{\operatorname{tors}}

\newcommand{\chr}[0]{\operatorname{char}}

\newcommand{\onto}[0]{\twoheadrightarrow}

\newcommand{\ann}[0]{\operatorname{Ann}}

\newcommand{\uhom}[0]{\operatorname{\underline{Hom}}}

\newcommand{\mor}[0]{\operatorname{Mor}}

\def\into{\DOTSB\lhook\joinrel\to}

\begin{document}
\bibliographystyle{amsalpha}

\title{Simultaneous normalization and algebra husks}
\author{J\'anos Koll\'ar}

\today

\maketitle

Let $X\to S$ be a morphism  with fibers
$\{X_s: s\in S\}$. We say that  $X\to S$ has a
simultaneous normalization if the 
normalizations of the fibers $\{\bar X_s: s\in S\}$
fit together to form a  flat family over $S$;
see (\ref{simnorm.defn}) for a precise definition.

The most famous result of this type, due to
 Hironaka \cite{hironaka}, says that
if $S$ is regular, the fibers $X_s$ are generically reduced and
the reductions of the   fibers $\red X_s$ are normal, then
 $\red X\to S$ is  flat  with normal fibers.
Several  related results are proved in  \cite{k-flat}.

For projective morphisms, a global condition
for the existence of the  simultaneous normalization
was developed by Chiang-Hsieh and Lipman \cite{ch-l}.
They consider the case when $S$ is normal with perfect residue fields
and $X\to S$ has reduced fibers, 
and prove that a  simultaneous normalization exists iff
the Hilbert polynomials of the normalizations of the fibers are all the same.
We give a proof in (\ref{ch-l.thm}).

When this condition fails,  a  simultaneous normalization exists
for some  subfamilies but not for others.
Our main result is  an analog of the 
Flattening decomposition theorem of  \cite[Lecture 8]{mumf66},
giving a precise description of
those subfamilies that have a simultaneous normalization:

\begin{thm}\label{main.thm}
 Let $f:X\to S$ be a  proper 
morphism whose fibers $X_s$ are  generically geometrically
reduced.

Then there is a morphism $\pi:S^n\to S$ such that
for any $g:T\to S$, the fiber product
$X\times_S T\to T$ has a 
simultaneous   normalization (\ref{simnorm.defn})
iff $g$ factors through $\pi:S^n\to S$.
\end{thm}

More precisely, we show that $\pi:S^n\to S$ represents
the functor of simultaneous   normalizations (\ref{simnorm.defn}).
We discuss in (\ref{simnorm.nonred.rem}) why the analogous result fails
if the fibers  are  not generically reduced.
The main result (\ref{main.P.thm}) establishes a similar
theorem for various partial normalizations.

The key technical step of the proof is to consider
 not just  the normalization  of a scheme $Y$,  but
all {\it algebra husks}  (\ref{husk.field.defn}) of $\o_Y$; equivalently, 
all finite birational maps $Y'\to Y$ as well.
More generally, algebra husks make sense for any coherent sheaf
of $\o_Y$-algebras  and they lead to  a well behaved
moduli functor (\ref{husk.exists.thm}).
  For arbitrary coherent sheaves this was
considered in \cite{hull-book}. It turns out to be easy to
derive the variant for  $\o_Y$-algebras from the results in
\cite{hull-book}.

In the last section 
we also settle two of the flatness conjectures from \cite{k-flat}
for proper morphisms.

\section{Algebra Husks}

We start by reviewing the notion of {\it husks} and their main
properties  \cite{hull-book}.

\begin{defn}\label{husk.field.defn}
  Let $X$ be a scheme over a field $k$ and $F$ a quasi coherent 
sheaf on $X$. Set $n:=\dim \supp F$.
A {\it husk} of $F$ is a  quasi coherent  sheaf $G$ together with a
homomorphism $q:F\to G$ such that
 $G$ has no associated primes of dimension $<n$, and
 $q$ is an isomorphism at all 
points of $X$ of dimension $=n$.

Such a  $G$ is also
an $\o_X/\ann(F)$ sheaf and so the particular choice of $X$
matters very little.

If, in addition, $F$ is a  quasi coherent  sheaf of $\o_X$-algebras, then a
 husk  $q:F\to G$ is called an {\it algebra husk}
if $G$ is a  quasi coherent sheaf of $\o_X$-algebras and 
$q$ is an algebra homomorphism.

Assume that $X$ is pure dimensional and generically reduced.
Then every coherent algebra husk of $\o_X$ is contained in
the structure sheaf $\o_{\bar X}$ of the normalization. 
If $X$ is of finite type over a field $k$ (more generally, 
if the local rings of $X$ are  {\it Nagata}
\cite[Sec.31]{mats-cr} or {\it universally Japanese}
 \cite[IV.7.7.2]{EGA}) then 
$\o_{\bar X}$ is coherent as an $\o_X$-module.
Thus, in these cases,  $\o_{\bar X}$ is 
the universal coherent algebra husk of $\o_X$.

If $X$ is not generically reduced, then there is no
universal coherent algebra husk. For instance,
for each $m\geq 0$,
$k[x,\epsilon]\to k\bigl[ x, x^{-m}\epsilon\bigr]$
is a coherent algebra husk.
\end{defn}



\begin{defn}\label{husk.defn}
  Let $f:X\to S$ be a morphism and $F$ a quasi coherent 
sheaf.
 Let $n$ be the relative dimension of
$\supp F\to S$.  
A {\it husk} of $F$ is a quasi coherent  sheaf $G$ together with a
homomorphism $q:F\to G$ such that
\begin{enumerate}
\item $G$ is flat over $S$ and
\item $q_s:F_s\to G_s$ is a husk   for every $s\in S$.
\end{enumerate} 


If, in addition, $F$ is a quasi coherent  sheaf of $\o_X$-algebras, then a
 husk  $q:F\to G$ is called an {\it algebra husk}
if $G$ is a quasi coherent  sheaf of $\o_X$-algebras and 
$q$ is an algebra homomorphism.

Note  that any multiplication map $m_G:G\otimes G\to G$
extending $m_F:F\otimes F\to F$ gives an algebra structure on $G$.
That is, $m_G$ is automatically associative
(and commutative if $m_F$ is). For example, associativity  is
equivalent to 
the vanishing of the difference map
$$
m_G(m_G,{\rm Id}_G)-m_G({\rm Id}_G,m_G):
G\otimes G\otimes G\to G.
$$
Since the target $G$ has no embedded points and the map vanishes
on a dense open set (since $m_F$ is associative),
the map is identically zero.

Note that husks and algebra husks are preserved by base change.
\end{defn}

\begin{defn} \label{husk.funct.defn} Let $f:X\to S$ be a morphism and 
$F$  a  coherent sheaf on $X$.
Let ${\it Husk}(F)(*)$ be the functor that to a
scheme $g:T\to S$ associates the
set of all  coherent husks of $g_X^*F$ with proper support over $T$,
where $g_X:T\times_SX\to X$ is the projection.

Assume that 
$H$ is an $f$-ample divisor and  $p(t)$  a polynomial. 
Let ${\it Husk}_p(F)(*)$ be the functor that to a
scheme $g:T\to S$ associates the
set of all  coherent husks of $g_X^*F$
with Hilbert polynomial $p(t)$.

If, in addition, $F$ is a coherent sheaf of $\o_X$-algebras, then
${\it Husk}^{\rm alg}(F)(*)$
(resp.\ ${\it Husk}_p^{\rm alg}(F)(*)$) denotes the functor that to a
scheme $g:T\to S$ associates the
set of all coherent  algebra husks of $g_X^*F$
(resp.\   coherent  algebra husks with Hilbert polynomial $p(t)$).
\end{defn}

\cite[Cor.12]{hull-book} shows that
${\it Husk}_p(F)$  has a fine moduli space ${\rm Husk}_p(F)$
which is an algebraic space over $S$.
Our basic  existence theorem  asserts that
similar results hold for algebra husks.

\begin{thm}\label{husk.exists.thm}
  Let $f:X\to S$ be a projective morphism, 
$H$ an $f$-ample divisor,  $p(t)$  a polynomial and 
 $F$  a  coherent sheaf  of $\o_X$-algebras.
Then 
\begin{enumerate}
\item ${\it Husk}_p^{\rm alg}(F)$ 
 has a fine moduli space ${\rm Husk}_p^{\rm alg}(F)$
which is an algebraic space of finite type over $S$.
\item The forgetful map 
 $\sigma: {\rm Husk}_p^{\rm alg}(F)\to {\rm Husk}_p(F)$ is a closed embedding.
\item If $F$ is flat at the generic points of $X_s\cap \supp F$
for every fiber $X_s$, then ${\rm Husk}_p^{\rm alg}(F)$
is proper over $S$.
\end{enumerate}
\end{thm}

Proof.  
For any $g:T\to S$, forgetting the algebra structure gives   a map
$$
\sigma_g:{\it Husk}^{\rm alg}_p(F)(T)\to {\it Husk}_p(F)(T).
$$
By  the first part of (\ref{alg.str.uniq.lem}), $\sigma_g$ 
is injective, that is,
a  husk admits at most one
structure as a sheaf of $\o_X$-algebras such that $F\to G$ is an algebra
homomorphism.
Next apply the second part of 
 (\ref{alg.str.uniq.lem})
to $S:={\rm Husk}_p(F)(T)$
to obtain that  ${\rm Husk}^{\rm alg}_p(F)(T):=S^{\rm alg}$
exists and is a closed subscheme of ${\rm Husk}_p(F)$,
proving (2). 

By  \cite[Cor.12]{hull-book} this in turn implies (1) and (3) since
 ${\rm Husk}_p(F)$ is  proper over $S$ 
if $F$ is flat at the generic points of $X_s\cap \supp F$
for every  $s\in S$.
\qed

\begin{lem}\label{alg.str.uniq.lem}
 Let $f:X\to S$ be a proper morphism,
 $F$ a coherent sheaf of $\o_X$-algebras
and $q:F\to G$ a coherent  husk. Then 
\begin{enumerate}
\item $G$ admits at most one
structure as an $\o_X$-algebra husk of $F$.
\item There is a closed subscheme $S^{\rm alg}\subset S$ such that
for a morphism $\pi:S'\to S$,
 $\pi_X^*G$ is an
$\o_{X'}$-algebra husk of  $\pi_X^*F$ iff $\pi$ factors through
$S^{\rm alg}$, where $\pi_X:S'\times_SX\to X$ is the projection.
\end{enumerate}
\end{lem}

Proof. We may assume that $S, S'$ are  affine.
The  first claim is also local on $X$; let
$U\subset X$ be affine.

The algebra structures are given by the multiplication maps
$m_F:F\otimes F\to F$ and $m_G:G\otimes G\to G$. 
There is an $h\in H^0(U,\o_U)$ which is
not a zero divisor on $G_s$ for every $s\in S$
 such that $h\cdot G|_U\subset F|_U$. 
For any  sections $g_1, g_2\in H^0(U, G|_U)$,
$h^2\cdot m_G(g_1\otimes g_2)=m_F(hg_1\otimes hg_2)$.
Since  multiplication by $h^2$ is
injective,  the above equality determines
$m_G(g_1\otimes g_2)$ uniquely, proving the first claim.

Next we prove (2)  in the 
projective case. 
For $r\gg 1$, there is an $h\in H^0(X,\o_X(r))$ that is
not a zero divisor on $G_s$ for every $s\in S$
 such that $h\cdot G\subset F(r)$. 
By the above considerations, the multiplication map
$m_F$ always extends to a multiplication map
$m_G^*:G\otimes G\to G(2r)$
and $G$ does  have an $\o_{X}$-algebra structure
iff $m_G^*$ actually maps $G\otimes G$ to $G\cong h^2\cdot  G
\subset  G(2r)$.

Note that since $G$ is flat over $S$ and $h^2$ is
not a zero divisor on $G_s$ for every $s\in S$,
 $G(2r)/G$ is also flat over $S$,
cf.\ \cite[Thm.22.5]{mats-cr}.
Hence we can choose $m>0$ such that 
$G(m)$ is generated by global sections and
$f_*\bigl(\bigl(G(2r)/G\bigr)(2m)\bigr)$
 is locally free and commutes with base change.

Pick generating sections
  $g_i\in H^0(X,G(m))$ and consider 
the composite
$$
g_i\otimes g_j\to m_G^*(g_i\otimes g_j)\to \bigl(G(2r)/G\bigr)(2m).
$$
By pushing forward, we obtain global sections
$$
\sigma_{ij}\in H^0\Bigl(S, f_*\bigl(\bigl(G(2r)/G\bigr)(2m)\bigr)\Bigr)
$$
such that 
for a morphism $\pi:S'\to S$,
 $\pi_X^*G$ is an
$\o_{X'}$-algebra husk of  $\pi_X^*F$ iff the $\pi^*\sigma_{ij}$
are all zero.
Thus 
$S^{\rm alg}:=\bigl(\sigma_{ij}=0\ \forall i,j)\subset S$
is the required subscheme.

In the proper but non-projective case, we 
first prove that
${\rm Husk}^{\rm alg}(F)$ exists and is of finite type over $S$.
To see this, we use the 
existence of the Hom-schemes of sheaves (\ref{hom.sch.exists}).

First consider $p:\uhom(G\otimes G, G)\to S$
with universal homomorphism
$u_p:p^*(G\otimes G)\to p^* G$. We can compose it with
$F\to G$ to obtain
$$
u'_p: p^*(F\otimes F)\to p^* G.
$$
Pulling back the multiplication map $m_F$ and the 
husk map $F\to G$ gives another homomorphism
$$
p^*m_F:p^*(F\otimes F)\to p^* F\to p^* G.
$$
We can view both of these maps as sections
$$
\uhom(G\otimes G, G)\rightrightarrows 
\uhom\bigl(p^*(F\otimes F), p^*G\bigr).
$$
Let $W\subset \uhom(G\otimes G, G)$ be the subscheme where
these two maps agree. Thus $W$ parametrizes
those husks for which the multiplication map
$F\otimes F\to F$ extends to a multiplication map
$G\otimes G\to G$. As noted in (\ref{husk.defn}), these are
the algebra husks of $F$.

Next we prove the  valuative criterion of   properness
for $p:W\to S$.
Let $T$ be the spectrum of a DVR with closed point $0\in T$
and generic point $t\in T$. Given $g:T\to S$
we have a husk $q_T:F_T\to G_T$ which is an algebra husk over
$t$. Set $Z:=\supp\coker q_0$. Then $F_0/\tors F_0\to G_0$ is an
isomorphism  over $X_0\setminus Z$, hence $G_T$ is an algebra
husk of $F_T$  over $X\setminus Z$.  

Note that
$G_0$ is $S_1$ and so $G_T$ is $S_2$ over its support
and $Z$ has codimension $\geq 2$ in $\supp G_T$.
In particular, every local section of $G_T$ over 
$X\setminus Z$ extends uniquely to a local section
over $X$ by \cite[III.3.5]{sga2}. Therefore, the multiplication map
$$
m_{X\setminus Z}:\bigl(G_T\otimes G_T\bigr)|_{X\setminus Z}
\to G_T|_{X\setminus Z}
$$
extends uniquely to a multiplication map
$m_X:G_T\otimes G_T\to G_T$.
Thus
$W\to S$  satisfies the 
 valuative criterion of   properness.

We have proved that   $p:W\to S$ is 
 a monomorphism of finite type that
satisfies  the valuative criterion of   properness.
 Thus $p:W\to S$ is a closed embedding and 
$S^{\rm alg}$ is its image.\qed

\begin{say}\label{hom.sch.exists}
(see \cite[III.7.7.8--9]{EGA}, \cite[4.6.2.1]{la-mb}, \cite[2.1.3]{lieb} or
\cite[33]{hull-book})
 Let $f:X\to S$ be proper. Let $M,N$ be coherent sheaves
on $X$ such  that $N$ is flat over $S$. Then there is a
separated   $S$-scheme of finite type 
$\uhom(M,N)$ parametrizing homomorphisms from $M$ to $N$.
That is, for any $g:T\to S$, there is a natural isomorphism
$$
\Hom_T(g_X^*M, g_X^*N)\cong \mor_S\bigl(T,  \uhom(M,N)\bigr),
$$ 
where $g_X: T\times_SX\to X$ is the fiber product
of $g$ with the identity of $X$.
\end{say}

\begin{exmp} This example shows that ${\rm Husk}^{\rm alg}(F)$
 is not always a union of connected components of  ${\rm Husk}(F)$.

 Consider the family of plane curves
$X:=(y^2z-x(x^2-tz^2))\subset \p^2\times \a^1\to \a^1$
over a field $k$.
Let $F:=\o_X$ and let $G$ be the subsheaf of 
rational functions generated by $\o_X$ and $y/x$. 
Since $x\cdot (y/x)=y$ and $y\cdot(y/x)=\frac1{z}\bigl(x^2-tz^2\bigr)$,
we see that  $G/F\cong k[t]$ and so $G$ is a husk of $F$.

Over the central fiber, $G_0$ is the coordinate ring
of the normalization of $X_0$, hence $G_0$ is an algebra husk of $F_0$.

For $t\neq 0$, $X_t$ is a smooth elliptic curve and
$G_t=\o_{X_t}(P_t)$ where $P_t$ denotes the origin.
The multiplication map gives a surjection
$$
G_t\otimes G_t\onto \o_{X_t}(2P_t)\supsetneq \o_{X_t}(P_t) = G_t.
$$
 Hence there is no algebra structure on $G_t$ extending the
algebra $\o_{X_t}$.
\end{exmp}

In the proper but non-projective case, we get the
following using \cite[Thm.39]{hull-book}.

\begin{thm}\label{husk.exists.prop.thm}
  Let $f:X\to S$ be a proper morphism and 
 $F$  a  coherent sheaf  of $\o_X$-algebras.
Then  ${\it Husk}^{\rm alg}(F)$ 
 has a fine moduli space ${\rm Husk}^{\rm alg}(F)$ 
and the forgetful map 
 ${\rm Husk}^{\rm alg}(F)\to {\rm Husk}(F)$ is a closed embedding.
\qed
\end{thm}

\begin{rem} As in \cite[Defn.9]{hull-book}, one can define the 
the functor ${\it QHusk}^{\rm alg}(F)$ of
algebra husks of quotients of $F$. It also has a fine
 moduli space ${\rm QHusk}^{\rm alg}(F)$
and the forgetful map 
 ${\rm QHusk}^{\rm alg}(F)\to {\rm QHusk}(F)$ is a closed embedding.
\end{rem}

\section{Simultaneous normalization}

\begin{defn}\label{simnorm.defn}
 Let $f:X\to S$ be a morphism.
A {\it simultaneous normalization} of $f$ is a morphism
$n:\bar X\to X$ such that
\begin{enumerate}
\item $n$ is finite and an isomorphism at the
generic points of the fibers of $f$, and
\item $\bar f:=f\circ n:\bar X\to S$ is flat with geometrically normal fibers.
\end{enumerate}

In characteristic 0, and over perfect fields, 
normal and geometrically normal are the same,
but over imperfect fields there are varieties which are
normal but not geometrically normal.

Note that, in general, a simultaneous normalization
need not be unique (\ref{norm.not.1}).

The functor of  simultaneous   normalizations
associates to a scheme $T\to S$ the set of 
simultaneous   normalizations of $X\times_ST\to T$.
\end{defn}

We start with a short proof of 
the  existence criterion \cite[Thm.4.2]{ch-l}.
The present form is somewhat more general since we
allow a semi-normal base and 
nonreduced fibers as well.
(For the definition of semi-normal, see \cite[I.7.2]{rc-book}.)

\begin{thm} \label{ch-l.thm}
Let $S$ be semi-normal with perfect residue fields at closed points.
 Let $f:X\to S$ be a projective 
morphism of pure relative dimension $n$ with generically 
reduced fibers. The following are equivalent:
\begin{enumerate}
\item $X$ has a simultaneous normalization
$n:\bar X\to X$.
\item The Hilbert polynomial of the normalization
of the fibers $\chi(\bar X_s,\o(tH))$
is locally constant on $S$.
\end{enumerate}
\end{thm}

Proof. The implication (1) $\Rightarrow$ (2) is clear.
To see the converse, we may assume that $S$ is connected.
Then $\chi(\bar X_s,\o(tH))$
is  constant; call it  $p(t)$. 
Set $S':=\red\bigl({\rm Husk}_p^{alg}(\o_X)\bigr)$
with universal family $X'\to S'$.
The structure sheaf $\o_X$ is flat over $S$
at the generic point of every fiber;
see \cite[II.2.3]{sga1} and \cite[Thm.8]{k-flat}
for the normal case and \cite[I.6.5]{rc-book} for the
semi-normal case.
Thus, by (\ref{husk.exists.thm}.3),
$\pi:S' \to S$ is proper.

Let $s\in S$ be a  closed point and
$\o_{X_s}\to F$  any algebra husk with
Hilbert polynomial $p(t)$. 
Since every coherent algebra husk of $\o_{X_s}$ is contained in the
structure sheaf of the normalization,
$$
p(t)=\chi(X_s,F(tH))\leq \chi(\bar X_s,\o_{\bar X_s}(tH))=p(t).
$$
Thus $F=\o_{\bar X_s}$ and since the residue field $k(s)$ is perfect,
this holds for any field extension of $k$. 
Therefore
$\pi$ is one-to-one  and surjective on closed geometric points.
Furthermore, the closed fibers of $X'\to S'$ are 
 geometrically normal, hence every fiber is
 geometrically normal.

Let now $g\in S$ be a generic point. By assumption,
$\chi(\bar X_g,\o_{\bar X_g}(tH))=p(t)$.
Thus $\o_{X_g}\to \o_{\bar X_g}$ is an algebra husk with 
Hilbert polynomial $p(t)$ and so 
the injection $g\into S$ lifts to $g\into S'$.
Therefore $S'\to S$ is an isomorphism
and $X'\to S'=S$  gives the simultaneous normalization.
\qed

\begin{exmp} The analog of (\ref{ch-l.thm}) fails for
semi-normalization, even for curves. 

As a simple example,
start with a flat family of curves $Y\to C$ whose general fiber
is smooth elliptic and $Y_c$ is a cuspidal rational curve for some $c\in C$.
Pick  2 smooth points in $Y_c$ and identify them to obtain
$X\to C$. The general fiber is still smooth elliptic but
$X_c$ has a cusp and a node plus an embedded point at the node.

The   semi-normalization of $X_c$ is a nodal rational curve, yet the
semi-normalizations do not form a flat family.
\end{exmp}

Next we state and prove our main result in a  general form.

\begin{defn}[Partial normalizations]\label{Pnorm.defn}
Let $P$ be a property of schemes or algebraic spaces
satisfying the following conditions.
\begin{enumerate}
\item $P$ is local, that is, $X$ satisfies $P$ iff an open cover satisfies $P$.
\item $P$ commutes with smooth morphisms, that is, if $X\to Y$ is smooth
and $Y$ satisfies $P$ then so does $X$.
\item If the  maximal dimensional generic points of $X$ satisfy $P$
then there is a unique smallest algebra husk $\o_X\to \bigl(\o_X\bigr)^P$
such that $X^P:=\spec_X \bigl(\o_X\bigr)^P$ satisfies $P$.
In this case $X^P\to X$ is called the {\it P-normalization} of $X$.
If $X^P=X$ then we say that $X$ is {\it P-normal}.
\item $P$-normalization is open. That is, 
given $X'\to X\to S$ such that  the composite $X'\to S$ is flat, 
 the set of points
$x\in X$ such that  $X'_s\to X_s$  is the $P$-normalization near $x$
 is open in $X$.
\end{enumerate}

Examples of such properties are:
\begin{enumerate}\setcounter{enumi}{4}
\item $S_1$, with $\o_X\to \o_X/(\mbox{torsion subsheaf})$ as the  
$P$-normalization.
\item $S_2$, with 
$\o_X\to j_*\bigl(\o_{X\setminus Z}/(\mbox{torsion subsheaf})\bigr)$
 as the  $P$-normalization  (or $S_2$-hull) where
$Z\subset X$ is a subscheme of codimension $\geq 2$ such that
$\o_{X\setminus Z}/(\mbox{torsion subsheaf})$ is $S_2$ and
$j:X\setminus Z\to X$ is the inclusion.
\item Normal, with the normalization.
\item Semi-normal, with the semi-normalization 
(cf.\ \cite[Sec.I.7.2]{rc-book}).
\item $S_2$ and semi-normal, with the $S_2$-hull of the semi-normalization.
\end{enumerate}

Note that in cases (5--6), generic points are always $P$-normal.
In cases (7--9) a generic point is  $P$-normal iff it is reduced. 

 Let $f:X\to S$ be a morphism.
As in (\ref{simnorm.defn}),
 a  {\it simultaneous P-normalization} of $f$ is a morphism
$n^P:(X/S)^P\to X$ such that
 $n^P$ is finite and an isomorphism at the
generic points of the fibers of $f$, the composite
 $f^P:=f\circ n^P:(X/S)^P\to S$ is flat 
and for every  geometric point $s\to S$
the induced map  $(X/S)^P_s\to X_s$ is the  $P$-normalization.

As before,
the functor of  simultaneous   $P$-normalizations
associates to a scheme $T\to S$ the set of 
simultaneous   $P$-normalizations of $X\times_ST\to T$.
\end{defn}

Our main technical theorem is the following.

\begin{thm}\label{main.P.thm}
 Let $f:X\to S$ be a  proper 
morphism.
Let $P$ be a property satisfying (\ref{Pnorm.defn}.1--4) and
assume that  the fibers $X_s$ are  generically geometrically
$P$-normal.

Then there is a morphism $\pi^P:S^P\to S$ 
that represents the functor of
simultaneous $P$-normalizations.

In particular, 
for any $g:T\to S$, the fiber product
$X\times_S T\to T$ has a 
simultaneous    $P$-normalization
iff $g$ factors through $\pi^P:S^P\to S$.

Furthermore,
$\pi:S^P\to S$ is one-to-one and onto on geometric points.
\end{thm}

Note that in cases (\ref{Pnorm.defn}.5--6), $f$ can be an arbitrary
proper morphism. In cases (\ref{Pnorm.defn}.7--9),
we assume that the fibers are generically geometrically reduced.
(The necessity of this condition is discussed in (\ref{simnorm.nonred.rem}).)
However, $f$ need not be flat nor  equidimensional.
\medskip

Proof.
By (\ref{husk.exists.prop.thm}), there is an algebraic space
${\rm Husk}^{\rm alg}(\o_X)$ parametrizing all algebra husks of
$\o_X$. Being a
$P$-normalization is an open condition, 
hence there is
an open subspace
$$
S^P:={\rm Husk}^{\rm P-n-alg}(\o_X)\subset {\rm Husk}^{\rm alg}(\o_X)
$$
parametrizing geometric $P$-normalizations  of $\o_X$.

 If $Y$ is an algebraic space over
an algebraically closed field then its $P$-normalization is unique
and is  geometrically $P$-normal.
This implies that $\pi$ is one-to-one and onto on geometric points. \qed

\begin{rem} In characteristic 0, this implies that  for every $s'\in S^P$,
$k(s')=k(\pi(s'))$, but in  positive characteristic
 $k(s')\supset k(\pi(s'))$ could be a purely inseparable
extension, even for the classical case of normalization.
For instance, if $k$ is a function field $K(t)$ of characteristic 3,
$S=\spec K(t)$ and $X$ the plane cubic
$(y^2z+x^3-tz^3=0)$ then $X$ is regular but not geometrically normal.
Over $K\bigl(\sqrt[3]{t}\bigr)$ it becomes singular
and its normalization is $\p^1$. 
Thus $S^n=\spec_{K(t)}K\bigl(\sqrt[3]{t}\bigr)$.
(If $P$=normal, we use $S^n$ to denote  $S^P$.)

As shown by (\ref{norm.not.1}),
$\pi$ need not be a locally closed embedding,
not even in characteristic 0.

Let $Y$ be a generically  geometrically reduced  scheme over
a field $k$ of positive characteristic.
In this case (\ref{main.P.thm}) implies that there is
a unique purely inseparable
extension $k'\supset k$ such that for any extension
$L\supset k$, the normalization  of $X_L$
is geometrically  normal iff $L\supset k'$.
\end{rem}

\begin{rem}\label{simnorm.nonred.rem} 
While the normalization of a nonreduced scheme
is well defined, it does not seem possible to define
simultaneous normalization for families with generically nonreduced fibers
over a  nonreduced base.

As  a simple example, let $\pi:\a^2\to \a^1$ be the projection to
the $x$-axis. Set $X:=(y^2-x^2=0)$. Then
$\pi:X\to \a^1$ has a nonreduced fiber over the origin.
The simultaneous normalization exists over $\a^1\setminus\{0\}$
and also over $(x=0)$ but not over any open neighborhood of $(x=0)$.
 What about over the nonreduced scheme $(x^n=0)$?
If we want to get a sensible functor, then there should not be a
simultaneous normalization over $(x^n=0)$ for large  $n$.

On the other hand, consider $Y_n:=\bigl(y^2-x^2=(y-x)^n=0\bigr)$.
This is the line $(y-x=0)$ with some embedded points at the origin.
The simultaneous normalization should clearly be the  line $(y-x=0)$.
If we want to get a functor, this should hold after base change to
 any subscheme
of $\a^1$.

Note, however, that $\pi:X\to \a^1$ and $\pi:Y_n\to \a^1$
are isomorphic to each other over $(x^n=0)$.
\end{rem}

\begin{exmp}[Simultaneous normalization not unique]\label{norm.not.1}

Even for flat families with reduced fibers, 
simultaneous normalization need not be  unique.
This, however, happens, only when the base is not reduced.

Let $k$ be a field  and consider
the trivial deformation  $k(t)[\epsilon]$.
If $D:k(t)\to k(t)$ is any derivation  then
$$
k[t]_D:=\bigl\{ f+\epsilon D(f): f\in k[t]\bigr\}+\epsilon k[t]
\subset k(t)[\epsilon]
$$
is a flat deformation of $k[t]$ over $k[\epsilon]$ which
agrees with the  trivial  deformation  iff $D(k[t])\subset k[t]$.
 
Consider the case $D(f(t)):=t^{-1}f'(t)$.
The deformation 
$k[t]_D$ is nontrivial since $D(t)=t^{-1}$.
On the other hand, 
$$
t^2=\bigl( t^2+\epsilon D(t^2)\bigr)-\epsilon\cdot  2,\qtq{and}
t^3=\bigl( t^3+\epsilon D(t^3)\bigr)-\epsilon\cdot  3t
$$
are both in $k[t]_D$, thus
$k[t]_D$ contains the trivial deformation
$$
k[t^2,t^3]+\epsilon k[t^2,t^3]
$$
of $k[t^2,t^3]$.
Hence both $k[t]_D$ and
$k[t]+\epsilon k[t]$ are simultaneous normalizations of
the trivial deformation of $k[t^2,t^3]$ over $k[\epsilon]/(\epsilon^2)$.

It is easy to see that $k[t]_D$
 cannot be extended to deformations over
$k[\epsilon]/(\epsilon^3)$, save  in characteristic 3,
where, for any $b\in k$,
$$
x+\epsilon\frac{1}{x} +\epsilon^2\Bigl(\frac{1}{x^3}+\frac{b}{x}\Bigr)
$$
generates an extension as a $k[\epsilon]/(\epsilon^3)$-algebra.
\end{exmp}

\begin{exmp}\label{simnorm.over.cusp} We give an example of surface
$f:X\to \spec k[t^2,t^3]$ such that
\begin{enumerate}
\item $X$ is reduced and $S_2$,
\item $f$ is flat except at a single point,
\item $\o_X$ has no hull,
\item over $\spec k[t]$, the hull is the structure sheaf of 
$\p^1\times \spec k[t]$ and
\item for $\chr k\neq 3$,  $\spec k[t]\to \spec k[t^2,t^3]$ represents the
simultaneous normalization functor.
\end{enumerate}

We start with the normalization of $X$, which is
$\p^1_{x:y}\times \a^1_t$. The map
$\p^1_{x:y}\times \a^1_t \to X$ will be a homeomorphism.
On the $y=1$ chart, $X$ is the spectrum of the ring
$$
R:=k\bigl[x^n+nx^{n-2}t: n\geq 2, x^nt^m: n\geq 0, m\geq 2\bigr]
\subset k\bigl[x,t\bigr].
$$
Note that $f(x)+g(x)t\in R$ iff $g(x)=\tfrac1{x}f'(x)$.

$R$ is finitely generated; one generating set is given by
$$
x^2+2t, x^3+3xt, t^2, xt^2, t^3, xt^3.
$$
Indeed, this set gives all the monomials
$t^m$ and $xt^m$ for $m\geq 2$. Now
$$
x^nt^m=x^{n-2}t^m\bigl(x^2+2t\bigr)-2x^{n-2}t^{m+1}
$$
gives all other monomials $x^nt^m$ for all $m\geq 2$.
Finally products of $x^2+2t$ and $ x^3+3xt$ give
 all the $x^n+nx^{n-2}t$ modulo $t^2$.

Consider $x^{-1}k[x]$ with the usual $k[x,t]$-module structure.
One easily checks that
$$
k[x,t]\to x^{-1}k[x]\qtq{given by}
f_0(x)+f_1(x)t+\cdots  \mapsto   f_1(x)-\tfrac1{x}f'_0(x)
$$
is an $R$-module homomorphism whose kernel is $R$.
Since $x^{-1}k[x]$ has no embedded points, we see that $R$ is $S_2$.

By explicit computation, the fiber of $R$ over the origin is
the cuspidal curve $k[x^2,x^3]$ with 2 embedded points at the origin.
Thus $R$ is generically flat over $k[t^2,t^3]$
but it is not flat at the origin.

The $(x=1)$ chart is easier.
It is given by the spectrum of the ring
$$
Q:=k\bigl[y^n-ny^{n+2}t: n\geq 1, y^nt^m: n\geq 0, m\geq 2\bigr]
\subset k\bigl[y,t\bigr].
$$
Note that $f(y)+g(y)t\in Q$ iff $g(y)=-y^3f'(y)$ and
$Q$ is flat over $ k[t^2,t^3]$.
\end{exmp} 

The next result shows that the 
above problems with simultaneous normalization
only appear in codimension 1 on the fibers.

\begin{prop}\label{normal.in.codim1.prop}
Let  $f:X\to S$ be  a proper and equidimensional morphism. 
Assume that there is   a closed subscheme $Z\subset X$ such that
\begin{enumerate}
\item  $\codim (X_s, Z\cap X_s)\geq 2$
 for every $s\in S$ and
\item  $X\setminus Z$ is flat over $S$ with geometrically normal fibers.
\end{enumerate}
Then  $\pi:S^n\to S$ as in (\ref{main.P.thm}) is a 
monomorphism. If $f$ is projective, then   $\pi:S^n\to S$ is a
locally closed decomposition,
 that is,
a locally closed embedding and a bijection on geometric points.
\end{prop}

Proof.  
First we show that for any $T\to S$, a simultaneous normalization
of $X_T:=X\times_ST\to T$ is unique. To see this, let
$h:Y_T\to X_T$ be a  simultaneous normalization
and $j:X_T\setminus Z_T\into X_T$ the open embedding.
Then $h_*\bigl(\o_{Y_T}\bigr)$ is an coherent sheaf on
$X_T$ which has depth 2 along $Z_T$ and
which agrees with $\o_{X_T}/(\mbox{torsion})$ 
outside  $Z_T$. Thus
$$
h_*\bigl(\o_{Y_T}\bigr)=j_*\bigl(\o_{X_T\setminus Z_T}\bigr),
\eqno{(\ref{normal.in.codim1.prop}.3)}
$$
which shows that $Y_T$ is unique.
Thus  $\pi:S^n\to S$  is a 
monomorphism.

In the projective case,
let $p(t)$ be the largest  polynomial that
occurs as a Hilbert polynomial of the normalization
of a geometric fiber of $f$ (\ref{norm.hp.finite.lem})
and let $S^n_p\subset S^n$ denote the open  subscheme
of normalizations with Hilbert polynomial $p$.
We prove that $S^n_p\to S$ is a proper monomorphism. 

Consider ${\rm Husk}^{\rm alg}_p(\o_X)\to S$. It parametrizes
partial normalizations of the fibers with Hilbert polynomial
$p(t)$. Since $p(t)$ is the  largest  Hilbert polynomial,
this implies that    ${\rm Husk}^{\rm alg}_p(\o_X)$
parametrizes normalizations of fibers that are geometrically normal. Thus
$S^n_p={\rm Husk}^{\rm alg}_p(\o_X)$ and so
$S^n_p\to S$ is  proper.

A proper monomorphism  is a closed embedding, 
hence $\pi:S^n_p\to S$ is a closed embedding.

Finally, we replace $S$ by $S\setminus \pi\bigl(S^n_p\bigr)$ and
conclude by Noetherian induction.\qed

\begin{lem} \label{norm.hp.finite.lem} (cf.\ \cite[Sec.3]{ch-l})
Let $f:X\to S$ be a projective morphism and $H$ an $f$-ample
divisor. For $s\in S$, let
$\chi(\bar X_{k(\bar s)},\o(tH))$
denote the  Hilbert polynomial of the normalization
of the geometric fiber of $f$ over $s$. Then
\begin{enumerate}
\item $s\mapsto \chi(\bar X_{k(\bar s)},\o(tH))$ is constructible.
\item If  $f$ has pure relative dimension $n$ with generically geometrically
reduced fibers then
$s\mapsto \chi(\bar X_{k(\bar s)},\o(tH))$ is upper semi continuous.
\end{enumerate}
\end{lem}

Proof. We may assume that $S$ is reduced.
Let $s\in S$ be a generic point, $K\supset k(s)$
an algebraic closure and $\bar X_K\to X_K$ the normalization.
There is a finite extension $L\supset k(s)$ such that
$\bar X_L\to X_L$ is geometrically  normal.
Let $S_L\to S$ be a quasi-finite morphism whose generic fiber
is $\spec L\to s$. Let $n_L:X_L^n\to \red\bigl(X\times_SS_L\bigr)$ be the
normalization. The generic fiber of $f\circ n_L$ is
geometrically normal. Thus, by shrinking $S_L$ if necessary,
we may assume that $f\circ n_L$ is flat with
geometrically normal fibers. In particular, the 
Hilbert polynomials of normalizations
of geometric fibers of $f$ are the same for every point of
the open set $\im(S_L\to S)$. The first part follows by Noetherian induction.

In order to prove  upper semi continuity,  it is enough to deal
with the case when $S$ is the spectrum of a DVR, $X$ is normal
and the generic fiber $X_g$ is geometrically normal.
Let $\bar X_{\bar 0}\to X_{\bar 0}$ denote the normalization of the 
geometric special fiber. Since $X$ is normal,
$X_0$ has no embedded points and the same holds for
$X_{\bar 0}$. We assumed that $X_0$ is generically geometrically
reduced, thus $X_{\bar 0}$ is generically reduced.
Thus $\o_{X_{\bar 0}}\to \o_{\bar X_{\bar 0}}$ is an injection
and 
$$
\chi(\bar X_{\bar 0}, \o(tH))\geq \chi(X_{\bar 0}, \o(tH))=
 \chi(X_{0}, \o(tH))= \chi(X_g, \o(tH)). \qed
$$

\section{Other Applications}

As another application, we prove  the flatness conjecture 
\cite[6.2.2]{k-flat} and a generalization of 
the conjecture
\cite[6.2.1]{k-flat}. The original conjectures are about
arbitrary morphisms, but here we have to restrict
ourselves to the proper case.
The example \cite[15.5]{k-flat} shows that
in (\ref{1st.flat.cor}) the rational singularity assumption is necessary.

\begin{cor}  \label{1st.flat.cor}
Let $S$ be a reduced scheme
 over a field of characteristic 0 and
 $f:X\to S$ a  proper  morphism. 
Assume that there is   a closed subscheme $Z\subset X$ such that
\begin{enumerate}
\item  $\codim (X_s, Z\cap X_s)\geq 2$
 for every $s\in S$,
\item  $X\setminus Z$ is flat over $S$ with normal fibers, and
\item the normalization $\bar X_s$ has rational singularities for
every $s\in S$.
\end{enumerate}
Let $j:X\setminus Z\into X$ be the injection and
$\bar X:=\spec_X j_*\bigl(\o_{X\setminus Z}\bigr)$.

Then $\bar f:\bar X\to S$ is flat and its  fibers are normal with only
  rational singularities.
\end{cor}

Proof.  The case when $S$ is the spectrum of a DVR is in
\cite[14.2]{k-flat}.

For $P=$normalization, let $\pi^n:S^n\to S$ be as in (\ref{main.P.thm}).
By (\ref{normal.in.codim1.prop}), $\pi^n$ is a monomorphism and by 
the above cited \cite[14.2]{k-flat}, 
$\pi^n$ satisfies the valuative criterion of properness,
hence it is proper.
Therefore,  $\pi^n$ is an isomorphism. The rest follows from
(\ref{normal.in.codim1.prop}.3).\qed

\begin{rem} The proof of (\ref{1st.flat.cor}) 
in fact shows that if a result of this type holds for a 
certain class of singularities (instead of rational ones)
when the base is the spectrum of a DVR,
then it also holds for an arbitrary reduced base.

In particular, it also applies
when the fibers have normal crossing singularities  in codimension one and
their $S_2$-hulls are semi-rational.
 The proof of \cite[14.2]{k-flat}
works in this case, using the semi-resolution theorem
of \cite{semi-res-book}.
\end{rem}

\begin{cor}  \label{2nd.flat.cor}
Let $0\in S$ be a normal, local scheme and 
 $f:X\to S$ a  proper morphism of relative dimension 1.
Assume that 
\begin{enumerate}
\item  $f$ is smooth at the generic points of $X_0$,
\item the generic fiber of $f$ is 
either smooth or defined over a field of characteristic 0   and
\item the reduced fiber $\red X_0$ has only finitely many
partial normalizations (\ref{f.m.part.norm}).
\end{enumerate}
Let $n:\bar X\to X$ denote the normalization.
Then $\bar f:=f\circ n:\bar X\to S$ is flat with reduced  fibers.
\end{cor}

Proof. 
Let $g\in S$ be the generic point. Then $\o_{X_g}$ is reduced, hence
 $\o_{\bar X_g}$ is a  husks of  $\o_{X_g}$ which gives an isolated 
point  $P_g\in {\rm Husk}^{\rm alg}(\o_{X_g})$.
Let $S'_g\subset {\rm Husk}^{\rm alg}(\o_{X_g})$ denote the
irreducible component containing $P_g$.
Then $S'_g$ is 0-dimensional
and  reduced  if $ \bar X_g=X_g$.
Thus its reduced  closure  
$S'\subset \red\bigl({\rm Husk}^{\rm alg}(\o_X)\bigr)$ 
is an irreducible component
such that
the induced map  $\pi: S'\to S$ is an isomorphism near $s$
if $X_g$ is smooth and a monomorphism in general.
Thus $\pi: S'\to S$ is birational
if the generic fiber of $f$ is 
either smooth or defined over a field of characteristic 0.

By (\ref{husk.exists.thm}.3), $\pi$ is proper.
The fiber $\pi^{-1}(0)$ parametrizes partial normalizations of
$X_0$, hence it is finite by assumption.
 Therefore, by Zariski's main theorem,
$\pi$ is an isomorphism. 

Let $u:Y'\to S'$ be the universal flat family of husks.
There is a natural morphism $Y'\to X$ which is finite and birational.

Since $S$ and the fibers of $u$ are $S_2$, so is
$Y'$ (cf.\ \cite[Thm.23.3]{mats-cr}).
  Moreover,  $u$ is smooth along the generic fiber
and along the generic points of the special fiber,
hence $Y'$ is regular in codimension 1. By Serre's criterion
$Y'$ is normal, hence $\bar X\cong Y$.
\qed

\begin{say}[Curves with finitely many
partial normalizations] \label{f.m.part.norm}

 We are interested in reduced curves
$C$ over a field $k$ such that only finitely many curves sit between
$C$ and its normalization, even after base change. 
That is, up to isomorphisms
which are the identity on $C$, 
there are only finitely many diagrams
$\bar C_{\bar k}\to C_i\to C_{\bar k}$. 
This condition depends on the singularities of $C$ only,
and there are only few singularities with this property.
By the results of \cite{gr-kn, ki-st, gr-kr},
the only such plane curve singularities
are the simple  singularities $A_n, D_n, E_6,E_7,E_8$.

Another series is given by the semi-normal curve
singularities. Over $\bar k$ these  are analytically
isomorphic to  the
coordinate axes in $\a^n$ for some $n$.

The example (\ref{quartic.exmp}) shows that ordinary quadruple points 
have infinitely many partial normalizations and
the conclusion of (\ref{2nd.flat.cor}) also fails for them.
\end{say}

\begin{say}(A correction to \cite[15.5]{k-flat}.)
Let $X$ be a smooth projective variety of dimension $n$
with $H^1(X,\o_X)\neq 0$ and $L$ an ample line bundle on $X$.
Let $C(X):=\spec \sum_m H^0(X, L^m)$ be the corresponding cone over $X$
with vertex $v\in C(X)$. 
Let $\pi:C(X)\to \a^n$ be a projection with 1-dimensional fibers.

The second part of  \cite[15.5]{k-flat} asserts that
$\pi$ is not flat at $v$ for $n\geq 2$. However, this holds only
when $n>2$. If $n=2$ then $\pi$ is flat but the fiber through $v$ has 
embedded points at $v$.

These also show that in (\ref{2nd.flat.cor})
some strong restrictions on the singularities are necessary.
\end{say}

\begin{exmp} \label{quartic.exmp}
Another interesting example is given by
the deformations of the plane quartic with an ordinary quadruple point 
$$
C_0:=\bigl(xy(x^2-y^2)=0\bigr)\subset \p^2.
$$
Let ${\mathbf C}_4\to \p^{14}$ be the universal family of
degree 4 plane curves and
${\mathbf C}_{4,1}\to S^{12}$ the 12-dimensional subfamily
whose general members are elliptic curves with 2 nodes.
$S^{12}$ is not normal, thus, to put ourselves in the settings of
(\ref{2nd.flat.cor}), we  normalize $S^{12}$ and
pull back the family.

We claim that if we take the normalization of the total space
$\bar{\mathbf C}_{4,1}\to \bar S^{12}$, we get a
 family of curves whose fiber over $[C_0]$ has embedded points.
Most likely, the family is not even flat, but I have not checked this.

We prove this by showing that in different families of curves
through $[C_0]\in S^{12}$ we get different flat limits.

To see this, note that the 
 semi normalization  $C_0^{sn}$ of $C_0$ can be thought of as
 4 general lines through a
point in $\p^4$. In suitable affine coordinates, its coordinate ring is
$$
k[u_1,\dots, u_4]/(u_iu_j: i\neq j)\supset k[u_1+u_3+u_4, u_2+u_3-u_4].
$$
There is a 1-parameter family of partial semi normalizations of $C_0$
corresponding to the 3-dimensional linear subspaces
$$
\langle u_1,\dots, u_4\rangle \supset W_{\lambda} \supset
\langle u_1+u_3+u_4, u_2+u_3-u_4 \rangle.
$$
Each $W_{\lambda}$ corresponds to  a projection of $C_0^{sn}$ to $\p^3$;
call the image $C_{\lambda}\subset \p^3$. Then
$C_{\lambda}$ is 4 general lines through a
point in $\p^3$; thus it is a $(2,2)$-complete intersection curve
of arithmetic genus 1.
(Note that the $C_{\lambda}$ are isomorphic to each other,
but the isomorphism will not commute with the map to $C_0$ in general.)
Every $C_{\lambda}$ can be realized as the special fiber
in a family $S_{\lambda}\to B_{\lambda}$ of $(2,2)$-complete intersection curves 
 in $\p^3$
whose general fiber is a smooth elliptic curve.

By projecting these families to $\p^2$, we get a 1-parameter family
$S'_{\lambda}\to B_{\lambda}$ of curves in $S^{12}$ whose special fiber is $C_0$. 

Let now $\bar{S}'_{\lambda}\subset  \bar{\mathbf C}_{4,1}$ be the preimage 
of this family in the normalization. Then $\bar{S}'_{\lambda}$
is dominated by the surface $S_{\lambda}$.

There are two possibilities. First, if $\bar{S}'_{\lambda}$
is isomorphic to $S_{\lambda}$, then the fiber of
 $\bar{\mathbf C}_{4,1}\to \bar S^{12}$ over $[C_0]$ is
$C_{\lambda}$. This, however, depends on $\lambda$, a contradiction.
Second, if $\bar{S}'_{\lambda}$ is not isomorphic to
 $S_{\lambda}$,
then the fiber of  $\bar{S}'_{\lambda}\to B_{\lambda}$ over the origin is
$C_0$ with some embedded points. Since $C_0$ has arithmetic genus 3,
we must have at least 2 embedded points.
\end{exmp}

 \begin{ack} I thank D.~Abramovich, G.-M.~Greuel, J.~Lipman and D.~Rydh 
for many useful comments and corrections.
Partial financial support  was provided by  the NSF under grant number 
DMS-0758275.
\end{ack}

\bibliography{refs}

\providecommand{\bysame}{\leavevmode\hbox to3em{\hrulefill}\thinspace}
\providecommand{\MR}{\relax\ifhmode\unskip\space\fi MR }
\providecommand{\MRhref}[2]{%
  \href{http://www.ams.org/mathscinet-getitem?mr=#1}{#2}
}
\providecommand{\href}[2]{#2}
\begin{thebibliography}{Mum66}

\bibitem[CHL06]{ch-l}
Hung-Jen Chiang-Hsieh and Joseph Lipman, \emph{A numerical criterion for
  simultaneous normalization}, Duke Math. J. \textbf{133} (2006), no.~2,
  347--390. \MR{MR2225697 (2007m:14003)}

\bibitem[GK85]{gr-kn}
G.-M. Greuel and H.~Kn{\"o}rrer, \emph{Einfache {K}urvensingularit\"aten und
  torsionsfreie {M}oduln}, Math. Ann. \textbf{270} (1985), no.~3, 417--425.
  \MR{MR774367 (86d:14025)}

\bibitem[GK90]{gr-kr}
G.-M. Greuel and H.~Kr{\"o}ning, \emph{Simple singularities in positive
  characteristic}, Math. Z. \textbf{203} (1990), no.~2, 339--354. \MR{MR1033443
  (90k:14001)}

\bibitem[Gro67]{EGA}
Alexander Grothendieck, \emph{\'{E}l\'ements de g\'eom\'etrie alg\'ebrique.
  {I--IV}.}, Inst. Hautes \'Etudes Sci. Publ. Math. (1960--67),
  no.~4,8,11,17,20,24,28,32.

\bibitem[Gro68]{sga2}
\bysame, \emph{Cohomologie locale des faisceaux coh\'erents et th\'eor\`emes de
  {L}efschetz locaux et globaux {$(SGA$} {$2)$}}, North-Holland Publishing Co.,
  Amsterdam, 1968, Augment\'e d'un expos\'e par Mich\`ele Raynaud, S\'eminaire
  de G\'eom\'etrie Alg\'ebrique du Bois-Marie, 1962, Advanced Studies in Pure
  Mathematics, Vol. 2. \MR{MR0476737 (57 \#16294)}

\bibitem[Gro71]{sga1}
\bysame, \emph{Rev\^etements \'etales et groupe fondamental.}, Lecture Notes in
  Mathematics, vol. 224, Springer Verlag, Heidelberg, 1971.

\bibitem[Hir58]{hironaka}
Heisuke Hironaka, \emph{A note on algebraic geometry over ground rings. {T}he
  invariance of {H}ilbert characteristic functions under the specialization
  process}, Illinois J. Math. \textbf{2} (1958), 355--366. \MR{MR0102519 (21
  \#1310)}

\bibitem[Kol95]{k-flat}
J{\'a}nos Koll{\'a}r, \emph{Flatness criteria}, J. Algebra \textbf{175} (1995),
  no.~2, 715--727. \MR{MR1339664 (96j:14010)}

\bibitem[Kol96]{rc-book}
\bysame, \emph{Rational curves on algebraic varieties}, Ergebnisse der
  Mathematik und ihrer Grenzgebiete. 3. Folge., vol.~32, Springer-Verlag,
  Berlin, 1996. \MR{MR1440180 (98c:14001)}

\bibitem[Kol08a]{hull-book}
\bysame, \emph{Hulls and husks},
  http://www.citebase.org/abstract?id=oai:arXiv.org:0805.0576, 2008.

\bibitem[Kol08b]{semi-res-book}
\bysame, \emph{Semi log resolution},
  http://www.citebase.org/abstract?id=oai:arXiv.org:0812.3592, 2008.

\bibitem[KS85]{ki-st}
K.~Kiyek and G.~Steinke, \emph{Einfache {K}urvensingularit\"aten in beliebiger
  {C}harakteristik}, Arch. Math. (Basel) \textbf{45} (1985), no.~6, 565--573.
  \MR{MR818299 (87d:14019)}

\bibitem[Lie06]{lieb}
Max Lieblich, \emph{Remarks on the stack of coherent algebras}, Int. Math. Res.
  Not. (2006), Art. ID 75273, 12. \MR{MR2233719 (2008c:14022)}

\bibitem[LMB00]{la-mb}
G{\'e}rard Laumon and Laurent Moret-Bailly, \emph{Champs alg\'ebriques},
  Ergebnisse der Mathematik und ihrer Grenzgebiete. 3. Folge., vol.~39,
  Springer-Verlag, Berlin, 2000. \MR{MR1771927 (2001f:14006)}

\bibitem[Mat86]{mats-cr}
Hideyuki Matsumura, \emph{Commutative ring theory}, Cambridge Studies in
  Advanced Mathematics, vol.~8, Cambridge University Press, Cambridge, 1986,
  Translated from the Japanese by M. Reid. \MR{MR879273 (88h:13001)}

\bibitem[Mum66]{mumf66}
David Mumford, \emph{Lectures on curves on an algebraic surface}, With a
  section by G. M. Bergman. Annals of Mathematics Studies, No. 59, Princeton
  University Press, Princeton, N.J., 1966. \MR{MR0209285 (35 \#187)}

\end{thebibliography}

\vskip1cm

\noindent Princeton University, Princeton NJ 08544-1000

\begin{verbatim}kollar@math.princeton.edu\end{verbatim}

\end{document}